\documentclass{amsart}

\usepackage{amsfonts}
\usepackage{amsmath}
\usepackage{amsthm}
\usepackage{amssymb}
\usepackage[english]{babel}
\usepackage{graphics}
\usepackage{latexsym}
\usepackage{longtable} \usepackage{mathrsfs}
\usepackage{graphicx}
\usepackage{wrapfig} 
\usepackage[pdftex,colorlinks=true]{hyperref}

\newtheorem{theorem}{Theorem}
\newtheorem{corollary}[theorem]{Corollary}
\newtheorem{lemma}[theorem]{Lemma}
\newtheorem{proposition}[theorem]{Proposition}

\newtheorem{claim}[theorem]{Claim}
\newtheorem{example}[theorem]{Example}

\theoremstyle{definition}
\newtheorem{definition}[theorem]{Definition}
\newtheorem{remark}[theorem]{Remark}

\newcommand{\mK}{\mathcal{K}}

\newcommand{\B}{\mathcal{B}}
\newcommand{\J}{\mathcal{J}}

\newcommand{\R}{\mathbb{R}}
\newcommand{\N}{\mathbb{N}}
\newcommand{\mS}{\mathbb{S}}
\newcommand{\mB}{\mathbb{B}}

\newcommand{\noi}{\noindent}
\newcommand{\ms}{\medskip}

\newcommand{\de}{\delta}
\newcommand{\De}{\Delta}
\newcommand{\e}{\varepsilon}
\newcommand{\si}{\sigma}

\newcommand{\Om}{\Omega}
\newcommand{\om}{\omega}

\newcommand{\Th}{\Theta}
\newcommand{\ze}{\zeta}

\newcommand{\av}{-\!\!\!\!\!\!\int}
\newcommand{\lharpoonup}{-\!\!\!\!\rightharpoonup}

\newcommand{\larrow}{\longrightarrow}
\newcommand{\ot}{\otimes}

\newcommand{\ri}{\rightarrow}

\newcommand{\p}{\partial}
\newcommand{\sub}{\subseteq}
\newcommand{\set}{\setminus}
\newcommand{\by}{\times}

\newcommand{\tr}{\textrm{tr}}

\newcommand{\diam}{\textrm{diam}}
\newcommand{\dist}{\textrm{dist}}

\newcommand{\Div}{\textrm{Div}}

\newcommand{\supp}{\textrm{supp}}

\newcommand{\bt}{\begin{theorem}}\newcommand{\et}{\end{theorem}}
\newcommand{\bd}{\begin{definition}}\newcommand{\ed}{\end{definition}}
\newcommand{\bl}{\begin{lemma}}\newcommand{\el}{\end{lemma}}
\newcommand{\beq}{\begin{equation}}\newcommand{\eeq}{\end{equation}}
\newcommand{\bc}{\begin{claim}}\newcommand{\ec}{\end{claim}}
\newcommand{\bex}{\begin{example}}\newcommand{\eex}{\end{example}}
\newcommand{\bcor}{\begin{corollary}}\newcommand{\ecor}{\end{corollary}}
\newcommand{\bp}{\begin{proof}}\newcommand{\ep}{\end{proof}}

\newcommand{\BPL}{\medskip \noindent \textbf{Proof of Lemma} }
\newcommand{\BPC}{\medskip \noindent \textbf{Proof of Claim} }
\newcommand{\BPCOR}{\medskip \noindent \textbf{Proof of Corollary} }
\newcommand{\BPP}{\medskip \noindent \textbf{Proof of Proposition} }

\numberwithin{equation}{section}

\begin{document}

\title[Nonsmooth Convex Functionals and Feeble Viscosity Solutions]{Nonsmooth Convex Functionals and Feeble Viscosity Solutions of singular Euler-Lagrange Equations}

\author{\textsl{Nikos Katzourakis}}
\address{Department of Mathematics and Statistics, University of Reading, Whiteknights, PO Box 220, Reading RG6 6AX, Berkshire, UK}
\email{n.katzourakis@reading.ac.uk}


\date{}



\begin{abstract} Let $F=F(A)$ be nonnegative, convex and in $C^2(\R^n \set \mK)$ with $\mK \subsetneqq \R^n$ a closed set. We prove that local minimisers in $(C^0\cap W^{1,1}_{\text{loc}})(\Om)$ of 
\[  \label{1} \tag{1}
E (u,\Om) \, :=\, \int_{\Om} F(Du),\ \  \Om \sub \R^n,
\]
are ``very weak" viscosity solutions on $\Om$ in the sense of Juutinen-Lindqvist \cite{JL} of the highly singular Euler-Lagrange equation of \eqref{1} expanded:
 \[  \label{2}
F_{AA}(Du):D^2u \, =  \, 0. \tag{2}
\]
The hypotheses on $F$ do not guarrantee existence of minimising weak solutions and include the singular $p$-Laplacian for $p\in(1,2)$. A much deeper converse is also true, if $\mK=\{0\}$ and extra natural assumptions are satisfied. Our main advance is that we introduce systematic ``flat" sup-convolution regularisations which apply to general singular nonlinear PDE in order to cancel the strong singularity of $F$. As an application we extend a classical theorem of Calculus of Variations regarding existence for the Dirichlet problem. These results extends previous work of Julin-Juutinen \cite{JJ} and Juutinen-Lindqvist-Manfredi \cite{JLM}.
\end{abstract}

\maketitle

\section{Introduction} \label{section1}

Let $F : \R^n \ri \R$ be a nonnegative convex function. Consider the functional $E :W^{1,1}_{\text{loc}}(\Om)\larrow [0,\infty]$ defined by
\beq \label{1.1}
E(u,\B)\,:=\, \int_\B F(Du),\  \text{ when\  }\ F(Du) \in L^1(\B), \ \B \sub \Om \text{ Borel},
 \eeq
with $E(u,\B):=\infty$ otherwise. In this paper we establish the equivalence between continuous local minimisers of $E$ in the space $W^{1,1}_{\text{loc}}(\Om)$ and of (appropriately defined) viscosity solutions in the sense of Crandall-Ishii-Lions \cite{CIL} of the Euler-Lagrange equation corresponding to \eqref{1.1} expanded:
\beq  \label{1.4}
F_{AA}(Du):D^2u =0.
\eeq
The notation is either self-explanatory or otherwise standard: $F_{A}$ and $F_{AA}$ stand for the 1st and 2nd derivatives of $F$, ``$:$" is the Euclidean inner product in $\R^{n \by n}$ and ``$\cdot$" is the inner product in $\R^n$. 
In this work a continuous local minimiser $u \in W^{1,1}_{\text{loc}}(\Om)$ is meant in the sense that 
\beq  \label{1.3}
E(u,\Om')\leq E(u+\psi,\Om'),\  \text{ for }\psi \in W^{1,1}_0(\Om'),\  \Om' \Subset \Om.
\eeq
In order to derive the PDE for local minimisers, we assume that
\beq \label{1.2}
\text{$F$ is convex, $F\geq F(0)=0$ and $F\in C^2(\R^n \set \mK)$, $\mK \subsetneqq \R^n$ closed.}  
\eeq
For the opposite direction, we will also need that
\beq \label{1.5}
\mK =\{F=0\}=\{0\},\ F \in C^1(\R^n)\ \text{ and } \, \int_0^1\!\!\av_0^t  \Big\{ \sup_{t<|a|<1}\De F(a)\Big\} dt\, ds < \infty
\eeq
and either that 
\beq  \label{1.6}
\text{the viscosity solution is locally Lipschitz continuous on }\Om,
\eeq
or that
\beq  \label{1.7}
\left\{
\begin{array}{c}
\big(F_A(b)-F_A(a)\big) \cdot (b-a)\, \geq \, c|b-a|^r,\\
\big| F_A(a)\big|\, \leq \, \dfrac{1}{c} |a|^{r-1},\\
\underset{|a|\ri \infty}{\lim\sup}\, \De F(a)\, <\, \infty,
\end{array}
\right.
\eeq
for some $c>0$, $r>1$ and $a,b \in \R^n$. We note that the standard example of $p$-Dirichlet density $F(A)=|A|^p$, $p>1$, satisfies all the above assumptions, even in the singular range of exponents $1<p<2$.

The primary advance in this paper compared to results already existing in the literature is that  \eqref{1.1} is a \emph{general nonsmooth functional} and \eqref{1.4} is \emph{highly singular}, since the Hessian of $F$ is undefined on $\mK$ and could even be unbounded near $\mK$. Under our assumption \eqref{1.2}, the Euler-Lagrange equation can not be treated in the classical setting of weak solutions. For, there is no way to infer that a function $u$ satsfying \eqref{1.3} is a distributional weak solution of
\beq   \label{1.12}
\Div\big(F_A(Du)\big) = 0
\eeq
since \eqref{1.1} may \emph{not} be Gateaux differentiable and $F_A(Du)$ may not even be measurable if $F \not\in C^1(\R^n)$. The appropriate extension of ``very weak" viscosity solutions which works in the case \eqref{1.4} has already been implicitely introduced in \cite{JLM} in the special case of the $p$-Laplacian and has been subsequently formalised in \cite{JL} as ``feeble viscosity solutions" (see also \cite{IS} by Ishii-Souganidis). This is nontrivial because if $F \not\in C^2(\R^n)$, extra caution is required since the PDE does not make sense on the set $\{Du\in \mK\}$. We introduce this definition, appropriate adapted to our case, later in Section \ref{section2}. Roughly, the idea is to view \eqref{1.4} as the \emph{free boundary} problem
\beq
\left\{
\begin{array}{l}
F_{AA}(Du):D^2u = 0, \ \text{ on } \Om \set \Om(u) , \ms\\
\ \ \ \ \ \ \ \ \ \ \ \ f(Du)=0, \  \text{ on }\Om(u),
\end{array}
\right.
\eeq
where $f=\dist(\cdot,\mK)$ and $\Om(u)= \{Du\in \mK\}$. This leads to ``ordinary" viscosity solutions on the set $\Om\set \Om(u)$, coupled by a 1st order differential inclusion ``$Du(x) \in \mK$" for $x\in \Om(u)$.

The results herein extend recent work of Julin-Juutinen \cite{JJ} done in the special case of the $p$-Laplacian, that is for $F(A)=|A|^p$ and respective equation
\beq  \label{1.10}
\De_p u \, =\, \Div\big(|Du|^{p-2}Du\big)=0,
\eeq
for $p>1$. This paper provides a simplified proof of the original result due to Juutinen-Lindqvist-Manfredi \cite{JLM} who proved the equivalence among three different notions of solution for the $p$-Laplacian, that of weak solutions based on integration-by-parts, that of viscosity solutions based on the maximum principle for the expanded version of the PDE
\beq  \label{1.11}
|Du|^{p-2}\left(I +(p-2)\frac{Du \ot Du}{|Du|^2}\right):D^2u\, =\, 0
\eeq
and an other based on nonlinear potential theory, introduced by Lindqvist in \cite{L} as an extension of the classical idea of Riesz for the Laplacian. Also, Juutinen-Lindqvist-Manfredi observed that the standard idea of semicontinuous envelopes does not work in the singular case when $1<p<2$ because we then allow for ``false" solutions. The primary advances of \cite{JJ} is that they bypass the heavy uniqueness machinery of Viscosity Solution theory which was employed in \cite{JLM}, and also consider the inhomogeneous case of $\De_p u=f$. 

Except for the more general setting that we consider in this paper, the main technical advance herein is that in the generality of \eqref{1.4} appropriate approximations have to be introduced in order to circumvent the strong singularities of $F_{AA}$. To this end, we introduce systematic regularisations of feeble viscosity solutions which we call ``flat sup/inf convolutions". The ``flat" counterparts   $(u^\e)_{\e>0}$ of the classical sup convolutions (see e.g.\ \cite{CIL,S}) are semiconvex approximations which satisfy additional estimates of the type
\beq \label{1.16}
D^2{u^\e} \geq -\frac{\Phi(|Du^\e|)}{\e}I
\eeq
for $\Phi$ is sufficiently ``flat" in order to cancel the singularity of $F_{AA}$ near $\mK=\{0\}$. The flatness property allows to show that if $u^\e$ is a strong subsolution of \eqref{1.4} a.e.\ on $\Om\set \{Du^\e=0\}$, then $u^\e$ is weak subsolution of \eqref{1.12} not only on $\Om\set \{Du^\e=0\}$  but on the whole of $\Om$. Roughly, this is achieved by choosing as $\Phi$
\beq  \label{phi}
\Phi(t) \, \approx \, \inf_{|a|>t} \frac{1}{\De F(a)}  ,\ \text{ for } t>0
\eeq
(for details see Sections \ref{section4} and  \ref{section5}). We believe that the tools developed in Section \ref{section4} are flexible and useful for nonlinear singular elliptic and parabolic equations in general, so this section is written independently of the PDE we are handling. Moreover, herein (as in \cite{JJ}) we do not appeal to the heavy uniqueness machinery of Viscosity Solutions. Our results and techniques extend to more general situations and to the case of nontrivial right-hand-side for \eqref{1.4}, but at the expense of added technical complexity. Hence, we decided to keep things as simple as possible, illustrating the main ideas.  

Let us now state the main result of this paper.

\begin{theorem} \label{th1} Let $F$ be a convex function on $\R^n$ satisfying \eqref{1.2}. Fix also $\Om\sub \R^n$ open. Then:

\ms

(a) Continuous local minimisers of \eqref{1.1} in $W^{1,1}_{\text{loc}}(\Om)$ are Feeble Viscosity Solutions of \eqref{1.4} on $\Om$ (Definition \ref{def2}).

\ms

(b) Conversely, if in addition $F$ satisfies \eqref{1.5} and also either \eqref{1.6} or \eqref{1.7},
then, Feeble Viscosity Solutions of \eqref{1.4}  on $\Om$ which are in $W^{1,1}_{\text{loc}}(\Om)$ (Definition \ref{def2}), are continuous weak solutions of \eqref{1.12}. Moreover, they are local minimisers of the functional \eqref{1.1} in $W^{1,r}_{\text{loc}}(\Om)$ (under \eqref{1.6} we have $r=1$ and under \eqref{1.7} we have $r>1$ as in the assumption).
\end{theorem}

We note that the conclusion of $(a)$ above remains true under the weaker assumption that $u$ is a \emph{spherical local minimiser in $W^{1,1}_{\text{loc}}(\Om)$}, namely when we use test functions supported on small balls and with small $W^{1,1}$ norm. As an application of $(a)$ above, in Corollary \ref{cor1} we obtain an extension of the classical theorem of Calculus of Variations regarding existence of solution to the Dirichlet problem
\beq   \label{1.14}
\left\{
\begin{array}{r}
F_{AA}(Du):D^2u =0, \ \ \text{ in }\Om,\ \, \ms\\
 u=b, \ \text{ on }\p \Om,
\end{array}
\right.
\eeq
where $\Om \Subset \R^n$, $b \in W^{1,1}_{\text{loc}}(\Om) \cap C^0(\overline{\Om})$ has finite energy on $\Om$ (that is $E(b,\Om)<\infty$) and $F$ satisfies \eqref{1.1} together with the strengthened coercivity
\beq   \label{1.15}
 F(A)\, \geq \, c|A|^s -\dfrac{1}{c},
\eeq
 for some $s>n$ and $c>0$. These assumptions are much weaker than those guaranteeing the existence of weak solutions. In the case of $(b)$, things are trickier and the extra assumptions \eqref{1.5}-\eqref{1.7} are required. They however provide the stronger conclusion that viscosity solutions of \eqref{1.4} are locally minimising weak solutions.

We conclude this introduction by noting some very interesting papers which relate to the results herein. In \cite{I}, Ishii considered the question of equivalence between weak and viscosity solutions of linear (degenerate elliptic) PDEs, in \cite{JLP} Juutinen-Lukkari-Parvianen consider the same question for the $p(x)$-Laplacian and in \cite{SV} Servadei and Valdinoci consider the same question for the fractional Laplacian. In \cite{JL}, Juutinen and Lindqvist consider the problem of the removability of level sets, which in the present case of our singular PDE \eqref{1.4} is relavant to the removability of $\{Du \in \mK\}$. Finally, we note that one further interesting non-smooth convex hamiltonian is
\[
F(A)\, =\, \max\big\{ |A|-1,0 \big\}^p, \ p>1,
\]
and relates to the problem of traffic congestion (\cite{SaV, CF}). This example was brought to our attention by one of the referees. Since the singular set $\mK$ here is a sphere (and not $\{0\}$), only $(a)$ of Theorem \ref{th1} and Corollary \ref{cor1} as they stand apply to this case. However, the convex hull of $\mK$ is the unit ball. Possible extensions of our results to this interesting case may be investigated in future work.

\section{Feeble Viscosity Solutions} \label{section2}

In this section we consider the appropriate adaptation of the definition of viscosity solutions for the singular PDE \eqref{1.4}. We will use as primary definition the version based on weak pointwise generalised derivatives (jets), rather than test functions. 

We begin by recalling from \cite{CIL} that for $u \in C^0(\Om)$, $\Om \sub \R^n$, the standard 2nd order subjets and superjets $\J^{2,\pm}u(x)$ of $u$ at $x\in \Om$ are defined as
\begin{align} \label{2.1}
\J^{2,+}u(x):= \Big\{(p,X) & \in \R^n  \! \by \mS(n)\ \Big|\ u(z+x)\leq u(x)+ p\cdot z \nonumber\\
&+\frac{1}{2}X  :z\ot z+o(|z|^2), \ \text{ as }z\ri0 \text{ in }\Om\Big\},
\end{align}
\begin{align}  \label{2.2}
\J^{2,-}u(x):= \Big\{(p,X) & \in \R^n  \! \by \mS(n)\ \Big|\ u(z+x)\geq u(x)+ p\cdot z \nonumber\\
&+\frac{1}{2}X  :z\ot z+o(|z|^2), \ \text{ as }z\ri0 \text{ in }\Om\Big\}
\end{align}
where $\mS(n):=\{A \in \R^{n \by n}: A_{ij}=A_{ji}\}$ denotes the symmetric $n\by n$ matrices.

\bd[Feeble Jets] \label{def1} Let $\mK \subsetneqq \R^n$  be a closed set. For $u \in C^0(\Om)$, $\Om \sub \R^n$, the \emph{2nd order Feeble Subjet $\J^{2,\pm}u(x)$ relative to $\mK$ of $u$ at $x\in \Om$} is defined as 
\begin{align} \label{2.3}
\J_\mK^{2,+}u(x):= \Big\{(p,X)\in \J^{2,+}u(x) \ \big|\ p \in \R^n \set \mK \Big\}
\end{align}
Similarly, the \emph{2nd order Feeble Superjet $\J^{2,\pm}u(x)$ relative to $\mK$ of $u$ at $x\in \Om$} is defined as 
\begin{align} \label{2.4}
\J_\mK^{2,-}u(x):= \Big\{(p,X)\in \J^{2,-}u(x) \ \big|\ p \in \R^n \set \mK \Big\}
\end{align}

\ed

\bd[Feeble Viscosity Solutions] \label{def2} Let $\mK \subsetneqq \R^n$ be a closed set and let
\beq \label{2.3a}
G\in C^0\big((\R^n \set \mK)\by \mS(n)\big).  
\eeq
Let also $u$ be in $C^0(\Om)$. We say that $u$ is a \emph{Feeble Viscosity Solution} of 
\beq  \label{2.5}
G(Du,D^2u) \geq 0
\eeq
on $\Om$ (or, subsolution of $G(Du,D^2u)=0$), when
\beq  \label{2.6}
\inf_{(p,X)\in \J_\mK^{2,+}u(x)}G(p,X) \geq 0,
\eeq
for all $x\in \Om$. Similarly, we say that $u$ is a \emph{Feeble Viscosity Solution} of 
\beq  \label{2.7}
G(Du,D^2u) \leq 0
\eeq
on $\Om$  (or, supersolution of $G(Du,D^2u)=0$), when
\beq \label{2.8}
\sup_{(p,X)\in \J_\mK^{2,-}u(x)} G(p,X) \leq 0,
\eeq
 for all $x\in \Om$. We say that $u$ is a \emph{Feeble Viscosity Solution of $G(Du,D^2u)=0$}  on $\Om$  when both \eqref{2.6} and \eqref{2.8} hold.
\ed

In the case of \eqref{1.4} we consider in this paper, we have $G(p,X)=F_{AA}(p):X$. We also note the obvious identity $\De F=F_{AA}:I=F_{A_i A_i}$.

\begin{remark} \label{rem2} For the sake of clarity, let us state also the equivalent point of view of viscosity solutions via touching test functions. We say that  $u$ is a Feeble Viscosity Solution of 
\eqref{2.5} on $\Om$, when for all $x\in \Om$ and $\psi \in C^2(\R^n)$ for which $u-\psi$ has a vanishing local maximum at $x$ and $D\psi(x)\not \in \mK$, we have
\beq \label{2.12}
G\big(D\psi(x),D^2\psi(x)\big) \geq 0.
\eeq
Similarly, we say that  $u$ is a Feeble Viscosity Solution of 
\eqref{2.7} on $\Om$, when for all $x\in \Om$ and $\psi \in C^2(\R^n)$ for which $u-\psi$ has a vanishing local minimum at $x$ and $D\psi(x)\not \in \mK$, we have
\beq \label{2.13}
G\big(D\psi(x),D^2\psi(x)\big)  \leq 0.
\eeq
 We say that $u$ is a \emph{Feeble Viscosity Solution} on $\Om$  when both \eqref{2.6} and \eqref{2.8} hold.
\end{remark}

\section{Local minimisers are Feeble Viscosity Solutions.} \label{section3}

In this section we establish the one half of Theorem \ref{th1}, packed in the following

\begin{proposition} \label{pr1} Fix $\Om\sub \R^n$ open and let $F$ be a convex function on $\R^n$ which satisfies \eqref{1.2} for some $\mK\subsetneqq \R^n$ closed. 

Then, continuous local minimisers of \eqref{1.1} in $W^{1,1}_{\text{loc}}(\Om)$ are Feeble Viscosity Solutions of \eqref{1.4}  on $\Om$ (Definition \ref{def2}).
\end{proposition}

As we have already mentioned, it suffices $u$ to be a \emph{spherical local minimiser} in the above result, but we will not use this generality.

\BPP \ref{pr1}. The argument utilised here follows an idea of Barron and Jensen from \cite{BJ}. Assume for the sake of contradiction that there is a $u \in(C^0\cap W^{1,1}_{\text{loc}})(\Om)$ which satisifes \eqref{1.3} but for some $x\in \Om$ \eqref{2.6} fails. Then, in view of Remark \ref{rem2} and standard arguments, there exists a smooth $\psi \in C^2(\R^n)$ and an $r>0$ such that $D\psi(x) \not\in \mK$ and
\beq   \label{3.1}
u-\psi < 0=(u-\psi)(x)
\eeq
on $\mB_r(x)\set \{x\}$, while there is a $c>0$ such that
\beq \label{3.2}
F_{AA}(D\psi(x)):D^2\psi(x) \leq -2c<0.
\eeq
In the standard way,  $\mB_r(x):=\{y\in \R^n:|y-x|<r\}$. Since $\R^n \set \mK$ is open and $\psi \in C^2(\R^n)$, we can decrease $r$ further to achieve
\beq
D\psi\big(\mB_r(x)\big) \, \sub \, \R^n \set \mK.
\eeq
Hence, the map $F(D\psi)$ is in $C^2\big(\mB_r(x)\big)$ and as such we have
\beq
F_{AA}(D\psi):D^2\psi\, =\, \Div\big(F_A(D\psi)\big)
\eeq
on $\mB_r(x)$. By restricting $r$ even further, \eqref{3.2} gives
\beq  \label{3.5}
-\Div\big(F_A(D\psi)\big)\, \geq c,
\eeq
on $\mB_r(x)$. By \eqref{3.1}, strictness of the maximum of $u-\psi$ implies that there is a $k>0$ small such that by sliding $\psi$ downwards to some $\psi-k$, we have
\beq \label{3.6}
\Om^+:=\big\{u-\psi+k>0 \big\} \sub \mB_r(x)\sub \Om
\eeq
and also $u=\psi-k$ on $\p \Om^+$. By multiplying \eqref{3.5} by $u-\psi+k \in W^{1,1}_0(\Om^+)$ and integrating by parts, we obtain
\beq \label{3.7}
\int_{\Om^+}F_A(D\psi) \cdot (Du-D\psi) \, \geq\, c\int_{\Om^+}|u-\psi+k|.
\eeq
Since $F$ is convex on $\R^n$, the elementary inequality
\beq \label{3.8}
F_A(a) \cdot (b-a) \, \leq\, F(b)-F(a)
\eeq
implies that \eqref{3.7} gives
\begin{align} \label{3.9}
c\int_{\Om^+}|u-\psi+k| \, &\leq\, \int_{\Om^+}F(Du) - \int_{\Om^+}F(D(\psi-k)) \nonumber\\
     &=\, E\big(u,\Om^+\big) - E\big(\psi-k,\Om^+\big).
\end{align}
In view of \eqref{1.3} and \eqref{3.6}, we have $E\big(u,\Om^+\big) - E\big(\psi-k,\Om^+\big)\leq 0$ and hence $\Om^+=\emptyset$, which is a contradiction. Hence, $u$ is a Feeble Viscosity Solution of 
\beq
F_{AA}(Du):D^2u\geq0
\eeq
on $\Om$. The supersolution property follows in the similar way and so does the proposition.
\qed

\ms

Proposition \ref{pr1} implies the following existence theorem:

\begin{corollary}[Existence for the Dirichlet Problem] \label{cor1} Assume that the convex function $F$ satisfies \eqref{1.2} together with \eqref{1.15} for some $s>n$, $\Om \sub \R^n$ is bounded and $b \in  W^{1,1}_{\text{loc}}(\Om)\cap C^0(\overline{\Om})$ has finite energy on $\Om$ i.e. $E(b,\Om)<\infty$ where $E$ is given by \eqref{1.1}. 

Then, the Dirichlet Problem \eqref{1.14} has a Feeble Viscosity Solution (Definition \ref{def2}), which is (globally) minimising for $E$ in $W^{1,s}_b(\Om)$.
\end{corollary}

\BPCOR \ref{cor1}. The argument is a simple implementation of the direct method of Calculus of Variations, which we include it for the sake of completeness. Since $E(b,\Om)<\infty$, by \eqref{1.15}  it follows that 
\beq \label{3.11}
C_2 \int_\Om |Db|^s \leq C_1|\Om|\, + \, \int_\Om F(Db)\leq C_1|\Om|+E(b,\Om)
\eeq
and hence by Poincar\'e inequality $b\in W^{1,s}(\Om)$. Thus, the infimum of $E$ in the affine space $W^{1,s}_b(\Om)$ is finite:
\beq \label{3.12}
0\, \leq\, e:= \inf_{W^{1,s}_b(\Om)}E\, \leq E(b,\Om)\, <\, \infty.
\eeq
Let $(u^m)_1^\infty$ be a minimising sequence. Since $E(u^m,\Om)\larrow e$ as $m\ri \infty$, by \eqref{1.15}, \eqref{3.11}, \eqref{3.12} and Poincar\'e inequality we have the uniform bound $\|u^m\|_{W^{1,s}(\Om)}\leq C$. Hence, there exists a subsequence along which we have $u^m \lharpoonup u$ as $m\ri \infty$ weakly in $W^{1,s}(\Om)$. In view of assumption \eqref{1.2}, the functional \eqref{1.1} is weakly lower-semicontinuous  in $W^{1,s}(\Om)$ (see e.g. Dacorogna \cite{D}, p. 94). Hence,
\beq
E(u,\Om)\, \leq \, \liminf_{m\ri \infty}E(u^m,\Om)\, =\,  e <\infty
\eeq
and as a result $u$ is minimiser of $E$ in $W_b^{1,s}(\Om)$. Since $s>n$, By Morrey estimate we have $u\in C^0(\Om)$. By Proposition \ref{pr1}, $u$ is a feeble viscosity solution of \eqref{1.4} and also $u=b$ on $\p \Om$. Hence $u$ solves \eqref{1.15} and the corollary follows.     \qed

\section{Flat Sup-Convolution Approximations.} \label{section4}

In this section we introduce the appropriate regularisations of viscosity solutions that allow to handle singular equations. For the reader's benefit, this section is independent of the rest of the paper since these regularisations are fairly general and might be useful in other contexts as well.

\begin{definition}[Flat sup/inf convolutions] \label{def3} Fix $\Th \in C^2(0,\infty)$ strictly increasing with 
\[
\Th(0^+)=\Th'(0^+)=\lim_{s\ri 0^+}\frac{\Th''(s)}{s}=0. 
\]
Given $u\in C^0(\overline{\Om})$, $\Om\Subset \R^n$ and $\e>0$, for $x\in \Om$ we define
\begin{align}
u^\e(x)\, :=\, \sup_{y\in \Om}\left\{u(y)-\frac{\Th(|x-y|^2)}{2\e}\right\}, \label{4.1}\\
u_\e(x)\, :=\, \inf_{y\in \Om}\left\{u(y)+\frac{\Th(|x-y|^2)}{2\e}\right\}.\label{4.2}
\end{align}
We call \emph{$u^\e$ the flat sup-convolution of $u$ and $u_\e$ the flat inf-convolution of $u$}.
\end{definition}

The nomenclature ``flat" owes to that the approximations satisfy flatness estimates of the type of \eqref{1.16} with $\Phi$ sufficiently flat at zero in order to cancel the singularity in the gradient variable of a nonlinear coefficient $G(Du,D^2u)$, when $u$ is a viscosity solution of such a PDE. The next result collects the main properties of $u^\e$ and $u_\e$.

\begin{lemma}[Basic Properties] \label{le1} Assume that 
\beq \label{4.3}
\text{$\Th \in C^2(0,\infty)$ is strictly increasing \&
$\Th(0^+)=\Th'(0^+)=\lim_{s\ri 0^+}\dfrac{\Th''(s)}{s}=0.$}
\eeq
If $u^\e$ and $u_\e$ are given by \eqref{4.1} and \eqref{4.2} respectively, for $\e>0$ we have:

\ms

(i) $u_\e=-(-u)^\e$ and $u^\e\geq u$ on $\Om$.

(ii) If we set 
\begin{align} 
u^\e_y(x) \, &:= \, u(y)-\frac{\Th(|x-y|^2)}{2\e}, \label{4.4}\\
X(\e)\, & :=\,  \Big\{ y \in \overline{\Om}\ \big|\ u^\e(x)=u^\e_y(x)\Big\} \label{4.5},
\end{align}
then $X(\e)$ is compact and for all $x^\e \in X(\e)$,
\beq \label{4.6}
|x-x^\e|\, \leq \, \sqrt{\Th^{-1}(4\|u\|_{C^0(\Om)}\e )}=:\rho(\e).
\eeq
That is, $X(\e)\sub \overline{\mB_{\rho(\e)}(x)}$. Moreover, 
\beq \label{4.7}
u^\e(x)\, =\, \max_{y\in \overline{\mB_{\rho(\e)}(x)}}\left\{u(y)-\frac{\Th(|x-y|^2)}{2\e}\right\}
\eeq

(iii) $u^\e \searrow u$ in $C^0(\overline{\Om})$ as $\e \ri 0$.

(iv) For each $\e>0$, $u^\e$ is semiconvex and hence twice differentiable a.e. on $\Om$. Moreover, if $(Du^\e, D^2u^\e)$ denote the pointwise derivatives, we have the estimate
\beq \label{4.8}
D^2u^\e \, \geq \, -\frac{1}{\e} \left(\sup_{0< t < d} \big\{2\Th''(t^2) t^2 + \Th'(t^2)\big\} \right) I
\eeq
a.e. on $\Om$, where $d:=\diam(\Om)$.

(v) We set
\beq \label{4.9}
\Om^\e\, :=\, \Big\{x\in \Om\ \big|\ \dist(x,\p \Om)>\rho(\e) \Big\}.
\eeq
If $u$ is a (Feeble) Viscosity Solution of 
\beq \label{4.10}
G(Du,D^2u)\geq 0
\eeq
on $\Om$ (Definition \ref{def2}), then $u^\e$ is a (Feeble) Viscosity Solution of 
\beq \label{4.11}
G(Du^\e,D^2u^\e)\geq 0
\eeq
on $\Om^\e$. Moreover, $u^\e$ is a strong solution, a.e.\ on $\Om^\e \set \{Du^\e \in \mK\}$.

(vi) (Magic properties) Assume in addition that
\beq \label{4.12}
\text{the function $T(t):= \Th'(t^2)t$ is strictly increasing on $(0,d)$, $d=\diam(\Om)$.}
\eeq
Then, 
\beq \label{4.13}
(p,X)\in \J^{2,+}u^\e(x) \ \ \Longrightarrow \ \ (p,X)\in \J^{2,+}u(x^\e),
\eeq
where
\beq \label{4.14}
X(\e) \ni x^\e = x +T^{-1}(\e p)
\eeq
and in \eqref{4.14} $T$ stands for the extension of $T$ on $\R^n$, that is $T(z):=\Th'(|z|^2)z$.

(vii) If  \eqref{4.12} holds, then for a.e.\ $x\in \Om$, the set $X(\e)$ of \eqref{4.5} is a singleton $\{x^\e\}$, and 
\beq   \label{4.15}
x^\e \, =\, x + T^{-1}\big(\e Du^\e(x)\big) .
\eeq

(viii)  If  \eqref{4.12} holds, then for a.e.\ $x\in \Om$, we have the estimate 
\beq  \label{4.16}
|Du^\e(x)|\, \geq \, \dfrac{1}{\e}T(|x-x^\e|),
\eeq 
where $x^\e$ is as in \eqref{4.15} and $T$  is as in \eqref{4.12}.

(ix)  If  \eqref{4.12} holds and $u\in W^{1,\infty}_{\text{loc}}(\Om)$, then for any $\Om' \Subset \Om$, we have 
\beq
\|Du^\e\|_{L^\infty(\Om'^\e)}\leq \|Du\|_{L^\infty(\Om')}.
\eeq
\end{lemma}

The proof is an extension of standard results in the literature for the ``ordinary" sup-convolutions correspoding to $\Th(t)=t$, but we provide it for the sake of completeness and for convenience of the reader. In particular, Lemma \ref{le1} extends results of \cite{JJ} in the special case of $\Th(t)=t^{{p}/{2(p-1)}}$ corresponding to the regularisation  method used for the singular $p$-Laplacian for $1<p<2$.

\BPL \ref{le1}. (i) is obvious. 

(ii) By (i), \eqref{4.1} and \eqref{4.4}, \eqref{4.5}, we have
\beq  \label{4.19}
u^\e(x) \, =\, u(x^\e)- \frac{\Th(|x^\e-x|^2)}{2\e}  \, \geq \, u(y) - \frac{\Th(|y-x|^2)}{2\e} ,
\eeq
for all $y\in \Om$. By choosing $y:=x$, we get 
\beq
\Th(|x^\e-x|^2)\, \leq\, 4\|u\|_{C^0(\Om)}\e.
\eeq

(iii) By assumption, $u$ is uniformly continuous on $\Om$. Hence, there is an increasing $\om \in C^0[0,\infty)$ with $\om(0)=0$ such that $|u(x)-u(y)|\leq \om(|x-y|)$, for all $x,y \in \Om$. Hence, by (ii) and \eqref{4.19}, for any $x\in \Om$,
\begin{align}
u^\e(x) \,= \, u(x^\e) \, \leq \, u(x)\, + \, \om(|x-x^\e|) \, \leq \, u(x)\, +\, \om\big(\rho(\e)\big),
\end{align}
while by (i) we have $u(x)\leq u^\e(x)$. In addition, it can be easily seen that $0<\e'<\e''$ implies $u^{\e''}\geq u^{\e'}$.

(iv) By \eqref{4.4}, we have
\begin{align} 
Du^\e_y(x)\, &=\, -\frac{\Th'(|x-y|^2)}{\e}(x-y), \label{4.21}\\
D^2u^\e_y(x)\, &=\, -\frac{1}{\e}\Big[\Th'(|x-y|^2)I + 2\Th''(|x-y|^2)(x-y)\ot (x-y) \Big]. \label{4.22}
\end{align}
By \eqref{4.4} and \eqref{4.7} we have $u^\e(x)=\sup_{y\in X(\e)} \{u^\e_y(x)\}$, while \eqref{4.22} readily implies
\begin{align}  \label{4.23}
D^2u^\e_y(x)\, &\geq \, -\frac{1}{\e}\left\{\sup_{y\in X(\e)}\Big\{\Th'(|x-y|^2) + 2\Th''(|x-y|^2)|x-y|^2 \Big\}\right\} I \nonumber\\
 &\geq \, -\frac{1}{\e}\left\{\sup_{0< t < d}\Big\{\Th'(t^2) + 2\Th''(t^2)t^2 \Big\}\right\} I,
\end{align}
for all $y\in X(\e)$. Hence, all $u^\e_y$ are semiconvex, uniformly in $y$. Since $u^\e$ is a supremum of semiconvex functions, the conclusion follows by invoking Alexandroff's theorem.

(v) By setting 
\beq   \label{4.24}
u^{\e,z}(x)\, :=\, u(x+z)-\frac{\Th(|z|^2)}{2\e}, 
\eeq
we immediately have that $\J^{2,+}u^{\e,z}(x)=\J^{2,+}u(x+z)$. Hence, for $|z|\leq \rho(\e)$, each $u^{\e,z}$ is a Feeble Viscosity Solution of 
\beq
G(Du^{\e,z} ,D^2u^{\e,z})\geq 0 
\eeq
on $\Om^\e$, if $u$ is a Feeble Viscosity Solution of $G(Du ,D^2u)\geq 0$ on $\Om$. Moreover, the classical result that pointwise suprema of Viscosity Subsolutions are Viscosity Subsolutions (\cite{CIL}, p.\ 23) extends to the Feeble case as well. For, it suffices to observe that by \eqref{4.24} and \eqref{4.1} we have 
\beq
u^\e(x)\, =\, \sup_{|z|< \rho(\e)}u^{\e,z}(x)
\eeq
and that since $\mK$ is closed, if $\psi$ touches $u^\e$ from above at $x \in \Om$ and $D\psi(x)\not\in \mK$, then there is a $\de>0$ such that $D\psi(\mB_\de(x)) \not\in \mK$.

(vi) Let $(p,X) \in \J^{2,+}u^\e(x)$. Then, there is a $\psi \in C^2(\R^n)$ such that $u^\e-\psi \leq (u^\e-\psi)(x)$ with $D\psi(x)=p$ and $D^2\psi(x)=X$. Hence, for all $z,y \in \Om$ and $x^\e \in X(\e)$,
\beq \label{4.25}
u(y)-\frac{\Th(|y-z|^2)}{2\e}-\psi(z)\, \leq\, u(x^\e)-\frac{\Th(|x^\e-x|^2)}{2\e}-\psi(x).
\eeq
For $y:=x^\e$,  we get
\beq
\Th(|z-x^\e|^2) + 2\e \psi(z)\, \geq\, \Th(|x-x^\e|^2) + 2\e \psi(x)
\eeq
and hence the function $z\mapsto \Th(|z-x^\e|^2) + 2\e \psi(z)$ has minimum at $x$ which implies that its gradient vanishes there. Consequencently,
\beq \label{4.27}
\Th'(|x^\e-x|^2)(x^\e-x)\,=\, \e p. 
\eeq
By \eqref{4.12}, we have that the operator $T:\R^n \ri \R^n$ given by $T(z):=\Th'(|z|^2)z$ is injective and hence \eqref{4.27} gives $x^\e  -x = T^{-1}(\e p)$, which is \eqref{4.14}. By \eqref{4.25} for $z:=y-x^\e+x$, we obtain
\beq
u(y)-\psi(y-x^\e+x)\, \leq\, u(x^\e)-\psi(x)
\eeq
which implies that $(p,X) \in \J^{2,+}u\big( x + T^{-1}(\e p) \big)$, as desired.         

(vii) Since $u^\e$ is semiconvex, for a.e.\ $x\in \Om$, we have $(Du^\e(x),D^2u^\e(x)) \in \J^{2,+}u^\e(x)$. The conclusion follows by  \eqref{4.14}.

(viii) Fix an $x\in \Om$ such that $Du^\e(x)$ exists. For any $e\in \R^n$ with $|e|=1$, we have
\begin{align} 
|Du^\e(x)|\, &\geq \, e \cdot Du^\e(x)  \nonumber\\
&=\, \frac{d}{dt} \Big|_{t=0}u^\e(x+te)  \nonumber\\
&= \, \frac{d}{dt}\Big|_{t=0} \left\{  \max_{y\in \overline{\mB_{\rho(\e)}(x)}}  \left( u(y)-\frac{\Th(|x+te-y|^2)}{2\e} \right)\right\}.  \nonumber
\end{align}
By Danskin's theorem (\cite{Da}) and  (vii), we obtain
\begin{align} \label{4.29}
|Du^\e(x)|\, &\geq \, \max_{y\in X(\e)}  \left\{ \frac{d}{dt}\Big|_{t=0} \left( u(y)-\frac{\Th(|x+te-y|^2)}{2\e} \right)\right\} \nonumber\\
&= \, \max_{y\in X(\e)}  \left\{ -\frac{1}{\e} \Th'(|x-y|^2) (x-y)\cdot e \right\}  \nonumber\\
&= \,  -\frac{1}{\e} \Th'(|x-x^\e|^2) (x-x^\e)\cdot e. \nonumber
\end{align}
If $x=x^\e$, \eqref{4.16} follows immediately since $T(0)=0$. If $x\neq x^\e$, we choose $e:=-\frac{x-x^\e}{|x-x^\e|}$ and then estimate \eqref{4.29} implies \eqref{4.16}.           

(ix) Let $p\in \J^{1,+}u^\e(x)$. Then, for a.e.\ $x\in \Om$, we have $p=Du^\e(x)$ and also $p\in \J^{1,+}u(x^\e)$, where $|x^\e-x|\leq\rho(\e)$. Fix $\Om'\Subset \Om$, $\e>0$ small and $x\in \Om'^\e$ and choose $\de>0$ small such that $\mB_\de(x^\e)\sub \Om'$. Since $u(x^\e+z)-u(x^\e)\leq p\cdot z+o(|z|)$ as $z\ri 0$, for $z:=-\de e$ with $|e|=1$, we have
\begin{align}
|p|\, =\, \max_{|e|=1}\, p \cdot e\, &\leq\, o(1)\, +\, \max_{|e|=1}\frac{u(x^\e) - u(x^\e-\de e)}{\de}  \nonumber\\
&\leq \, o(1)\, +\, \sup_{x,y \in \mB_\de(x^\e)}\frac{|u(y)-u(x)|}{|y-x|}\\
&\leq \, o(1)\, +\, \|Du\|_{L^\infty(\Om')}, \nonumber
\end{align}
as $\de \ri 0$. Hence, for a.e.\ $x\in \Om'^\e$, we have $|Du^\e(x)|\leq  \|Du\|_{L^\infty(\Om')}$. \qed

\ms

The following result contains the main new property of the approximations \eqref{4.1}, \eqref{4.2}, which fails for the standard sup/inf convolutions and justifies the necessity of a general $\Th$ function.

\begin{lemma}[Flatness Estimates]  \label{le2}  Let $\Phi \in C^0(0,\infty)$ be a strictly increasing function with $\Phi(0^+)=0$ such that
\beq \label{Osgood}
\int_0^1 \dfrac{dt}{\Phi(t)}\, <\, \infty.
\eeq
Fix a domain $\Om \Subset \R^n$ and set $d:=\diam(\Om)$.

Then, there exists a strictly increasing function $\Th \in C^2(0,\infty)$ satisfying
\beq \label{Theta}
\Th(0^+)=\Th'(0^+)=\lim_{s\ri 0^+}\frac{\Th''(s)}{s}=0
\eeq
and such that, if we set $T(t):=\Th'(t^2)t$, the function $T'$ is positive and increasing on $(0,d)$. 

Moreover, the sup-convolution operator given for any $u\in C^0(\overline{\Om})$ by \eqref{4.1} satisfies the properties (i)-(ix) of Lemma \ref{le1} together with the estimate
\beq \label{4.30}
D^2u^\e \, \geq \, -\frac{\Phi (|Du^\e |) }{\e}  I
\eeq
a.e.\ on $\Om$, for any $\e>0$.
\end{lemma}

\BPL \ref{le2}. \textbf{Step 1.}  We define a function $T \in C^1(0,\infty)$ as follows: consider the initial value problem
\beq \label{4.33}
\left\{
\begin{array}{l}
T'(t)\, =\, \Phi\big( T(t)\big),\ \  t>0,\ms\\
\ T(0)\, =\, 0.
\end{array}
\right.
\eeq
In view of our assumption \eqref{Osgood}, Osgood's non-uniqueness criterion of ODE theory implies that \eqref{4.33} has a \emph{nontrivial} solution $T$, which is positive and strictly increasing on $(0,\infty)$, since $\Phi>0$ on $(0,\infty)$. Moreover, the composition $T'=\Phi \circ T$ is strictly increasing as well and $T'(0^+)=\Phi(0^+)=0$.

\ms

\noi \textbf{Step 2.}  We define $\Th \in C^2(0,\infty)$ by
\beq \label{4.34}
\Th(t)\, :=\, 2\int_0^{\sqrt{t}} T(s)ds.
\eeq
Obviously, $\Th(0^+)=0$ and since $T'\geq T'(0^+)=0$, we get 
\beq
\Th(t) \, \leq\, 2\sqrt{t} \, T(\sqrt{t})\, = \, \sqrt{t} \, o(\sqrt{t})\, =\, o(t), 
\eeq
as $t\ri 0$, which implies $\Th'(0^+)=0$. By differentiating \eqref{4.34}, we have the identity
\beq  \label{4.36}
\Th'(t^2)t\,=\, T(t)
\eeq
which implies $\Th>0$ and $\Th'>0$  on $(0,\infty)$. Moreover, the identity
\beq \label{4.39}
2\Th''(t^2)t^2\, +\, \Th'(t^2) \, =\, T'(t) 
\eeq
implies $\lim_{s\ri 0^+}\Th''(s)/s=0$ and hence \eqref{Theta} ensues.
\ms 

\noi \textbf{Step 3.} Let $\Th$ be defined by \eqref{4.34}, \eqref{4.33}. Then, $\Th$ as well as $T$ (given by  \eqref{4.36}) satisfy all the assumption of Lemma \ref{le1}. Hence, the sup-convolution operator defined by \eqref{4.1} for this $\Th$ satisfies the properties (i)-(ix) of Lemma \ref{le1}.

\ms

\noi \textbf{Step 4.} We now establish \eqref{4.30}. Fix $u\in C^0(\overline{\Om})$ for an $\Om\Subset \R^n$ and $0<\e<1$. For a.e.\ $x\in \Om$, $u^\e$ is twice differentiable at $x$. Fix such an $x\in \Om$. Then, by \eqref{4.16}, we have
\beq \label{4.38}
|x-x^\e| \, \leq\, T^{-1}\big(\e |Du^\e(x)| \big).
\eeq
On the other hand,  \eqref{4.39} and \eqref{4.8} imply
\begin{align} \label{4.40}
D^2u^\e(x)\, &\geq \,-\frac{1}{\e}\Big( 2\Th''(|x-x^\e|^2)|x-x^\e|^2+\Th'(|x-x^\e|^2)  \Big) I  \nonumber\\
                     &=\, -\frac{1}{\e}T'(|x-x^\e|) I.
\end{align}
Since $T'$ is increasing, by \eqref{4.40} and \eqref{4.38} we obtain 
\beq \label{4.41}
D^2u^\e \, \geq \, -\frac{1}{\e} T'\left(  T^{-1}\big(\e |Du^\e |\big) \right) I,
\eeq
a.e.\ on $\Om$. Finally, since $T'\circ T^{-1}$ is increasing, by \eqref{4.41}, \eqref{4.33} and by using that $\e |Du^\e |\leq |Du^\e |$, we have
\begin{align}
D^2u^\e \, & \geq \, -\frac{1}{\e} T'\left(  T^{-1}\big(|Du^\e |\big) \right) I \nonumber\\
&  = \, -\frac{1}{\e} \Phi\Big(T\left(  T^{-1}\big(|Du^\e |\big) \right) \Big) I\\
& = \, -\frac{1}{\e} \Phi\big(|Du^\e |\big)   I,\nonumber
\end{align}
a.e.\ on $\Om$. The lemma ensues.         \qed \ms

\begin{remark} Our assumption \eqref{1.5} on the radial integrability the Laplacian $\De F$ near the origin will guarrantee that the function $\Phi$ (roughly given by \eqref{phi}) has the ``Osgood property" \eqref{Osgood} which is needed for the construction of the flat supconvolutions.
\end{remark}

\section{Feeble viscosity solutions are weak locally minimising solutions.} \label{section5}

In this section we utilise the systematic approximations of Section \ref{section4} to establish the second half of Theorem  \ref{th1}. The first half has been established in Proposition \ref{pr1}.
 
\ms

\noi \textbf{Motivation of the method.}  Roughly, the idea of the usage of flat sup-convolutions is the following: choose the function 
\[
\Phi(t) \, = \, \inf_{|a|>t} \frac{1}{\De F(a)}, 
\]
and consider the flat sup-convolution $u^\e$ for the respective $\Th$. If $u^\e$ is a strong subsolution a.e.\ on $\Om\set \{Du^\e=0\}$ of 
\[
F_{AA}(Du^\e):D^2u^\e \, \geq\, 0,
\]
the flatness estimate
\[
D^2u^\e \, \geq \, -\frac{\Phi (|Du^\e |) }{\e}  I
\]
gives the lower bound
\[
F_{AA}(Du^\e):D^2u^\e \, \geq\ -\frac{\Phi (|Du^\e |) }{\e} F_{AA}(Du^\e):I\, \geq \, -\frac{1}{\e}
\]
and application of Fatou lemma allows to infer
\[
- \int_\Om F_A(Du^\e) \cdot D\psi\, \geq\, \int_{\Om\set\{Du^\e=0\}}  \psi F_{AA}(Du^\e):D^2u^\e\, \geq\, 0,
\]
for non-negative test functions. However, the above reasoning is too simplistic to apply exactly as it stands. Several regularisations are required in order to make this work, and this causes substantial  complications. For the simpler case of the $p$-Laplacian, see \cite{JJ}.

\ms

The main result here is

\begin{proposition} \label{pr2} Let $\mK\subsetneqq \R^n$ be closed and $F$ a convex function on $\R^n$ satisfying \eqref{1.2}, \eqref{1.5} and also either \eqref{1.6} or \eqref{1.7}. Fix also $\Om\sub \R^n$ open.

Then, if $u \in W^{1,1}_{\text{loc}}(\Om)$ is a Feeble Viscosity Solution of \eqref{1.4} on $\Om$ (Definition \ref{def2}), then $u$ is a continuous weak solution of \eqref{1.12} and also local minimiser of the functional \eqref{1.1} in $W^{1,r}_{\text{loc}}(\Om)$. Moreover, under \eqref{1.6} we have $r=1$ and under \eqref{1.7} we have $r>1$ as in the assumption.
\end{proposition}

The following lemma is the first step towards the proof of Proposition  \ref{pr2}.

\begin{lemma} \label{le3} Let $F : \R^n \ri \R$ be convex and satisfy assumptions \eqref{1.2} and \eqref{1.5}. 
\ms

(i) There exists a positive strictly decreasing function $\rho\in C^0(0,\infty)$ with $\rho(0^+)=\infty$ and $\rho(\infty)>0$ which depends only on $F$ such that, for each $R>1$, the function $\Phi^R \in C^0(0,R)$ given by
\beq \label{5.16}
\Phi^R(t)\, :=\, \underset{t<|a|<R}{\inf} \Big\{\frac{1}{\rho(|a|)+\De F(a)}\Big\} ,
\eeq
is positive, strictly increasing and satisfies $\Phi^R(0^+)=0$ and
\beq \label{Phi^R}
\int_0^1\dfrac{dt}{\Phi^R(t)}\, <\, \infty.
\eeq 
If in addition $F$ satisifes
\[
\underset{|a|\ri \infty}{\lim\sup}\, \De F(a)\, <\, \infty,
\]
then we may take $R=\infty$ and $\Phi^\infty \in C^0(0,\infty)$ is also positive, strictly increasing and satisfies the same properties.
\ms

(ii) Fix a domain $\Om\sub \R^n$. Assume that $v\in C^0(\Om)$ is a semiconvex strong subsolution of
\beq \label{5.1}
F_{AA}(Dv):D^2v\, \geq\, 0
\eeq
a.e.\ on $\Om\set \{Dv=0\}$. Moreover, suppose that $\|Dv\|_{L^\infty(\Om)}< R$ and that for some $\e>0$, $v$ satisfes the flatness estimate
\beq  \label{5.3}
D^2v\, \geq \, -\frac{\Phi^R(|Dv|)}{\e}I
\eeq
a.e.\ on $\Om$. Then, it follows that $v$ is a weak subsolution of 
\beq  \label{5.2}
\Div\big(F_A(Dv) \big)\,\geq\,0
\eeq
on $\Om$. 

If $\|Dv\|_{L^\infty(\Om)}=\infty$, then the same conclusion follows if in addition we have ${\lim\sup}_{|a|\ri \infty}\, \De F(a)< \infty$ and \eqref{5.3} is satisfied for $\Phi^\infty$. 
\end{lemma}

The function $\rho$ is explicitely constructed in the proof and is a correction term which arises because we need to mollify $F$ near the singularity at $\{0\}$. The above lemma has a symmetric counterpart for supersolutions which we refrain from stating explicitely. 

\BPL \ref{le3}. \textit{Proof of $(i)$:}

\textbf{Step 1.} We begin by utilising the assumption \eqref{1.5} in order to construct an explicit modulus of differentiability for $F$ in terms of $\De F$. Note that by assumption $F\in C^1(\R^n)\cap C^2(\R^n\set \{0\})$, zero is a strict global minimum for $F$ and $\De F\geq 0$ on $\R^n\set \{0\}$. We set
\beq \label{omega}
\om(t)\, :=\, \sqrt{t}\, +\, \int_0^t\Big\{\sup_{s<|a|<1}\De F(a) \Big\}ds.
\eeq
The ``$\sqrt{t}$" term is needed in order to avoid small technical complications which arise if $\De F$ is not unbounded and strictly radially decreasing near the origin.
\ms

\noi \textbf{Claim.} \textit{The function $\om$ given by \eqref{omega} is strictly increasing, concave and in $C^1(0,1)$. Moreover, $\om(0^+)=0$ and also satisfies
\beq \label{omegaprop}
t \mapsto \frac{\om(t)}{t} \text{ is strictly decreasing, }\int_0^1\frac{\om(t)}{t}dt < \infty.
\eeq
Morevoer, for any $a \in \R^n$ with $0<|a|<1$, we have
\beq  \label{5.6}
0  \leq F(a)\, \leq\, \om(|a|)|a|,\ \ \ \big|F_A(a) \big| \leq C\om(|a|),
\eeq
where $C=C(n)>0$ depeneds only on the dimension.}
\ms

\noi \textit{Proof of Claim.} By \eqref{1.5}, there is a $C>0$ such that
\begin{align}
\infty > \, C \, \geq & \  \int^1_{\frac{1}{2}} \av_0^t \Big\{\sup_{s<|a|<1}\De F(a) \Big\}ds \, dt \nonumber\\
\geq & \ \frac{1}{2} \av_0^{t_0} \Big\{\sup_{s<|a|<1}\De F(a) \Big\}ds \nonumber
\end{align}
for some $t_0 \in [1/2, \, 1]$. Hence, we have the estimate
\[
\int_0^{\frac{1}{2}} \Big\{\sup_{s<|a|<1}\De F(a) \Big\}ds\, \leq \, 2t_0C\, \leq 2C
\]
which implies that the positive strictly decreasing function 
\[
\de(s)\, :=\, \frac{1}{2\sqrt{s}}\, +\, \sup_{s<|a|<1}\De F(a) 
\]
is in $(L^1\cap C^0)(0,1)$. As a result, $\om$ is concave, strictly increasing, $\om(0^+)=0$ and $\om'(t)=\de(t)>0$ for $0<t<1$. Moreover,
\[
\left(\frac{\om(t)}{t}\right)'\, =\, \frac{1}{t^2}\left[t\, \de(t)\, -\, \int_0^t\de(s)ds \right]\, <\, 0,
\]
for $0<t<1$, and also (by assumption \eqref{1.5})
\[
\int_0^1\frac{\om(t)}{t}dt\, =\, \int^1_0\frac{1}{t}\int_0^t \left\{\frac{1}{2\sqrt{s}}\, +\, \sup_{s<|a|<1}\De F(a) \right\}ds \, dt\, <\, \infty.
\]
Thus, \eqref{omegaprop} has been established. Finally, by Taylor's theorem, for any $a \in \R^n$ with $0<|a|<1$ and $0<\e<1$,
\begin{align}
\Big|F(a)-F(\e a)-(1-\e)F_A(\e a)\cdot a \Big|\, = & \, \left| \int_\e^1 (1-t)F_{AA}(ta):a\ot a \, dt \right|\nonumber\\
   						\leq & \, |a|^2\int_\e^1 (1-t)\Big|F_{AA}(ta):\frac{a\ot a}{|a|^2} \Big|dt \nonumber\\
 						\leq & \, |a|^2\int_\e^1 (1-t)\De F(ta) dt. \nonumber
\end{align}
Hence,
\begin{align}
\Big|F(a)-F(\e a)-(1-\e)F_A(\e a)\cdot a \Big|\, \leq & \ |a|^2\int_\e^1 \Big\{\sup_{t<s<1}\De F(sa)\Big\} dt \nonumber\\
						= & \ |a|\int_{\e|a|}^{|a|} \Big\{\sup_{t<|A|<|a|}\De F(A)\Big\} dt \nonumber\\
					\leq & \ |a|\int_0^{|a|} \Big\{\sup_{t<|A|<1}\De F(A) \Big\} dt, \nonumber
\end{align}
which gives
\[
\Big|F(a)-F(\e a)-(1-\e)F_A(\e a)\cdot a \Big|\, \leq  \, |a|\, \om(|a|),
\]
for $0<|a|<1$ and $0<\e<1$. By passing to the limit as $\e \ri 0$, we obtain the desired estimate
$0  \leq F(a)\leq \om(|a|)|a|$. Similarly, for any $e\in \R^n$, $|e|=1$, we have
\beq \label{est1}
\big|\big(F_A(a)-F_A(\e a)\big)\cdot e \big|\, \leq \, |a|\int_\e^1 \Big|F_{AA}(ta):\frac{a}{|a|}\ot e\Big| \, dt.
\eeq
By norm equivalence on $\mS(\R^n) \sub \R^{n \by n}$, there is a $C=C(n)>0$ such that, for any non-negative symmetric $n\by n$ matrix $X$,
\begin{align} \label{est2} 
\Big|X:\frac{a}{|a|}\ot e \Big|\, &\leq \, \max_{|E|=1}\big\{X:E\big\} \nonumber\\ 
&=\,  |X| \\
& \leq \, C(n) \, \max_{|e|=1}\big\{X:e\ot e\big\} \nonumber\\
& \leq \, C(n)\, \tr(X). \nonumber
\end{align}
In view of the estimates \eqref{est1} and \eqref{est2} and by arguing as before, we conclude that 
\[
\big|F_A(a)-F_A(\e a)\big|\, \leq  \, C(n)\, \om(|a|),
\]
for $0<|a|<1$ and $0<\e<1$. Hence by letting $\e \ri 0$ we see that \eqref{5.6} has been established and the proof of the claim is complete.  \qed
\ms

\noi \textbf{Step 2.} We begin by introducing appropriate $C^2$ approximations of $F$ in $C^1(\R^n)$. Fix $0<\de<1$ and choose $\ze \in C^\infty[0,\infty)$ such that $\ze \equiv 0$ on $[0,1/2]$ and $\ze \equiv 1$ on $[1,\infty)$. Set
\beq \label{5.4}
\ze^\de(a)\, :=\, \ze \Big( \frac{|a|}{\de}\Big), \ \ a\in \R^n.
\eeq
Then, $\ze^\de \in C^\infty(\R^n)$, $\ze^\de \equiv 0$ on $\mB_{\de/2}(0)$ and $\ze^\de \equiv 1$ on $\R^n \set \mB_{\de}(0)$ and
\beq \label{5.5}
\big|D\ze^\de \big|\, \leq \, \frac{C}{\de}\chi_{\mB_{\de}(0)},\ \ \ \big|D^2\ze^\de \big|\, \leq \, \frac{C}{\de^2}\chi_{\mB_{\de}(0)},
\eeq
for some universal $C>0$. We set $F^\de:=\ze^\de F$. Then, we have that $F^\de \in C^2(\R^n)$, $F^\de\equiv 0$ on $\mB_{\de/2}(0)$ and $F^\de \equiv F$ on $\R^n \set \mB_{\de}(0)$. For $|a|\leq \de$, we have
\beq \label{5.7}
\big|F(a)-F^\de(a) \big|\, =\, \big|(1-\ze^\de(a)) F(a)\big|\, \leq\, \om(|a|)|a|\, \leq\, \om( \de) \de,
\eeq
and also
\begin{align} \label{5.8}
\big|F_A(a)-F^\de_A (a) \big|\, &\leq \, |D \ze^\de(a)| |F(a)| \, +\, | F_A(a)| |1-\ze^\de(a)| \nonumber\\
&\leq\, \frac{C}{\de}\om(|a|)|a|\, +\, C\om(|a|)\\ 
&\leq \, C\om(\de),\nonumber
\end{align}
as $\de\ri 0$. Hence, $F^\de \ri F$ in $C^1(\R^n)$ as $\de\ri 0$, as desired.

\ms

\noi \textbf{Step 3.} Now we construct the function $\rho$ of the statement.  By \eqref{5.5} and \eqref{5.6}, for $0<|a|<1$ and $0<\de<1$ we have
\begin{align} \label{5.9}
\big| \De F^\de(a) \big|\, &=\, \big|F(a)\De \ze^\de(a)\, +\, 2 D\ze^\de(a) \cdot F_A(a)\, +\, \ze^\de(a)\De F(a)\big| \nonumber\\
&\leq \om(|a|)|a|\, \frac{C}{\de^2}\chi_{\mB_\de(0)}(a)\, +\, 2\frac{C}{\de}\chi_{\mB_\de(0)}(a)\, \om(|a|)\, +\, \De F(a)\\
& \leq\, C\frac{\om(\de)}{\de}\chi_{\mB_\de(0)}(a)\, +\, \De F(a) \nonumber\\
& \leq\, C\frac{\om(|a|)}{|a|}\, +\, \De F(a). \nonumber
\end{align}
The last inequality owes to that $t\mapsto \om(t)/t$ is strictly decreasing on $(0,1)$. We set
\beq \label{5.10}
\rho(t)\,:=\, C\frac{\om(t)}{t} ,\ \ 0<t<1,
\eeq
and extend it on $[1,\infty)$ as a strictly decreasing positive function (for example, set $\rho(t):=\frac{C}{2}\om(1)(e^{1-t}+1)$ for $t\geq 1$). By \eqref{5.9} and \eqref{5.10} we have
\beq \label{5.11}
\sup_{0<\de<1}\big|\De F^\de(a) \big| \, \leq \, \rho(|a|)\, +\, \De F(a),
\eeq
for $a\neq 0$. Fnally, we employ \eqref{5.16} to define for $R>1$ the function $\Phi^R \in C^0(0,R)$ which is positive and strictly increasing with $\Phi^R(0^+)=0$. By utilising \eqref{omegaprop}, we obtain the estimate
\begin{align}
\int_0^1\frac{dt}{\Phi^R(t)}\, &=\, \int_0^1\left\{ \sup_{t<|a|<R}\Big[\rho(|a|)\, +\, \De F(a)\Big]\right\} dt \nonumber\\
&\leq\, \int_0^1  \Big\{\sup_{t<|a|<R}\rho(|a|)\Big\} dt \, +\, \int_0^1  \Big\{\sup_{t<|a|<R} \De F(a)\Big\} dt \nonumber\\
& \leq\, C\int_0^1 \frac{\om(t)}{t}dt \, +\,  \om(1) \,+\, \sup_{1<|a|<R}\De F(a)\nonumber\\
& <\, \infty. \nonumber
\end{align}
Hence, \eqref{Phi^R} has been established as well. If ${\lim\sup}_{|a|\ri \infty}\, \De F(a)< \infty$, then for $R=\infty$ the function $\Phi^\infty \in C^0(0,\infty)$ is also positive, increasing with  $\Phi^\infty(0^+)=0$ and satisfies the same estimate.

\ms

\noi \textit{Proof of $(ii)$:}

\noi \textbf{Step 1.} We now establish that the semiconvex function $v$ satisfies the inequality
\beq  \label{5.19}
-\int_\Om D\psi \cdot F^\de_A(Dv)\, \geq\, \int_\Om \psi F^\de_{AA}(Dv):D^2v,
\eeq
for all $\psi \in C^\infty_c(\Om)$, $\psi\geq 0$ and $0<\de<1$. We note that this is not a trivial application of integration by parts due to the existence of the singular part of the full Hessian measure of $v$ and equality in \eqref{5.19} may fail. Since $v$ is semiconvex on $\Om$, $(Dv,D^2v)$ exist  a.e.\ on $\Om$ and also there is $C>0$ such that $D^2v \geq -\frac{C}{\e}I$ a.e.\ on $\Om$. To prove \eqref{5.19}, we regularise $v$ further in the standard way by convolution, that is for $\si>0$ we consider the mollifier $v*\eta^\si$, which is also semiconvex uniformly in $\si>0$. 

Since by semiconvexity we have that $v \in W^{1,\infty}_{\text{loc}}(\Om)$ (e.g.\ \cite{EG}, p.236), by Dominated Convergence it follows that
\beq
D\psi \cdot F^\de_A(Dv*\eta^\si) \, \larrow \, D\psi \cdot F^\de_A(Dv), \label{5.19a}
\eeq
 in $L^1(\Om)$, as $\si \ri0$. Moreover,
\beq \label{5.19b}
F^\de_{AA} (Dv*\eta^\si):(D^2v*\eta^\si)\,  \larrow \, F^\de_{AA} (Dv):D^2v,
\eeq
a.e.\ on $\Om$ as $\si \ri 0$, and also for $\si>0$ small we have the $L^1$ lower bound
\beq \label{5.19c}
\psi F^\de_{AA} (Dv*\eta^\si):D^2v*\eta^\si \, \geq \, -\frac{C}{\e} \psi \max_{|a|\leq \|Dv\|_{L^\infty(\supp(\psi))}}|\De F^\de(a)|.
\eeq
Hence, by \eqref{5.19a}, \eqref{5.19b}, \eqref{5.19c} and Fatou Lemma, we have
\begin{align}
-\int_\Om D\psi \cdot F^\de_A(Dv)\, &=\, -\lim_{\si \ri 0} \int_\Om D\psi \cdot F^\de_A(Dv*\eta^\si)  \nonumber\\
&=\, \underset{\si\ri 0}{\lim\inf} \int_\Om \psi F^\de_{AA} (Dv*\eta^\si):(D^2v*\eta^\si)\\
&\geq\, \int_\Om \underset{\si\ri 0}{\lim\inf}\, \psi F^\de_{AA} (Dv*\eta^\si):(D^2v*\eta^\si)  \nonumber\\
&=\,  \int_\Om \psi F^\de_{AA}(Dv):D^2v. \nonumber
\end{align}
Hence, \eqref{5.19} follows.

\ms

\noi \textbf{Step 2.} Since by $(i)$ we have $F^\de\ri F$ in $C^1(\R^n)$ and  also $F^\de \equiv 0$ on $\mB_{\de/2}(0)$, we may use Dominated Convergence theorem to pass in the limit as $\de \ri 0$ in  \eqref{5.19}:
\begin{align} \label{5.24}
-\int_\Om D\psi \cdot F_A(Dv)\, &=\, -\underset{\de\ri 0}{\lim} \int_\Om D\psi \cdot F^\de_A(Dv) \nonumber\\
&\geq\, \underset{\si\ri 0}{\lim\inf}  \int_\Om \psi F^\de_{AA} (Dv):D^2v  \\
&=\, \underset{\si\ri 0}{\lim\inf}  \int_{\Om \set \{Dv=0\}} \psi F^\de_{AA} (Dv):D^2v. \nonumber
\end{align}
We now utilise the assumption \eqref{5.3} in order to pass to the limit as $\de\ri 0$. Let $\Phi^R$ be given by \eqref{5.16}. The choice of $R>0$ depends on our assumptions on $F$ and $v$:
\beq \label{5.17}
\left\{
\begin{array}{l}
R:=\|Dv\|_{L^\infty(\Om)} +1,\ \text{\ \  if \ } Dv \in L^\infty(\Om),\ms\ms\\
R:=\infty, \text{\hspace{63pt}  if \ } \underset{|a|\ri \infty}{\lim\sup}\, \De F(a) < \infty. 
\end{array}
\right.
\eeq
Then, by \eqref{5.11} we have
\begin{align} \label{5.25}
F^\de_{AA} (Dv):D^2v \, &\geq\,-\frac{1}{\e}\Phi^R(|Dv|) F^\de_{AA} (Dv):I \nonumber\\
&=\,-\frac{1}{\e}\Phi^R(|Dv|) \De F^\de (Dv) \\
&\geq\,-\frac{1}{\e}\Phi^R(|Dv|) \Big( \rho(|Dv|)+ \De F (Dv)\Big) \nonumber\\
&\geq\,-\frac{1}{\e}, \nonumber
\end{align}
a.e.\ on $\Om$, and this gives an $L^1$ lower bound in order to use Fatou Lemma in \eqref{5.24}. Hence, we have
\begin{align} \label{5.26}
-\int_\Om D\psi \cdot F_A(Dv)\, &\geq\,  \int_{\Om \set \{Dv=0\}} \underset{\de\ri 0}{\lim\inf} \, \psi F^\de_{AA} (Dv):D^2v\\
&\geq\,  \int_{\Om \set \{Dv=0\}}  \psi F_{AA} (Dv):D^2v. \nonumber
\end{align}
By our assumption on $v$, the right hand side of \eqref{5.26} vanishes. Hence,
\beq
-\int_\Om D\psi \cdot F_A(Dv)\, \geq\, 0
\eeq
for any $\psi \in C^\infty_c(\Om)$, $\psi\geq 0$. The lemma follows.        \qed \ms

\begin{remark}
The $L^1$ lower bound estimate \eqref{5.25} is the main reason for the necessity to introduce the sup-convolution approximations in Section \ref{section4} which satisfy flatness properties that cancel the singularity of the coefficients. See also \cite{JJ} where this idea has already been implicitely utilised in the special case of the singular $p$-Laplacian when $1<p<2$ and $F(A)=|A|^p$.  
\end{remark}

We may now prove the proposition.

\BPP \ref{pr2}. Let $u\in (C^0 \cap W^{1,1}_{\text{loc}})(\Om)$ be a Feeble Viscosity Solution of $F_{AA}(Du):D^2u\geq 0$ on $\Om\sub \R^n$. Fix $\phi \in C^\infty_c(\Om)$ with $\phi \leq 0$  and $\Om'  \Subset \Om$ such that $\supp(\phi)\sub \Om'$. Obviously, $u \in C^0(\overline{\Om'})$.

\ms

\noi \textbf{1st case under assumption \eqref{1.6}.} Since $u\in W^{1,\infty}_{\text{loc}}(\Om)$, we have $\|Du\|_{L^\infty(\Om')}<\infty$. Consider the increasing function $\Phi^R \in C^0(0,\infty)$ with $\Phi^R(0^+)=0$ defined by \eqref{5.16} in Lemma \ref{le3}, where as $R$ we take $\|Du\|_{L^\infty(\Om')} +1$. We then extend $\Phi^R$ on $(0,\infty)$ by setting 
\beq
\Phi(t)\, :=\, 
\left\{
\begin{array}{l}
\Phi^R(t),\ \ \ \ \ \ 0<t<R,\ms\\
\Phi^R(R)\dfrac{t}{R},\ \ t \geq R.
\end{array}
\right.
\eeq
Consider for $\e>0$ the flat sup-convolution $u^\e$ of $u$ restricted on $\Om'$, as given by \eqref{4.1}, where as $\Th \in C^2(0,\infty)$ we take the function given by Lemma \ref{le2} for the selected $\Phi$.  Then $u^\e$ satisfies all the properties of Lemmas \ref{le1}, \ref{le2}. Hence the flatness estimate \eqref{4.30} holds. Since $u\in W^{1,\infty}_{\text{loc}}(\Om)$, we have
\beq \label{5.28}
\|Du^\e\|_{L^\infty(\Om'^\e)}\, \leq\, \|Du\|_{L^\infty(\Om')}
\eeq
($\Om'^\e$ as in Lemma \ref{le1}). Moreover, by Lemma \ref{le1}, $u^\e$ is a strong subsolution of
\beq
F_{AA}(Du^\e):D^2u^\e \geq 0,
\eeq
a.e.\ on $\Om'^\e \set \{Du^\e=0\}$. Since $u^\e$ is semiconvex, it satisfies the assumptions of Lemma   \ref{le3}. Hence, $u^\e$ is a locally Lipschitz weak subsolution on $\Om'^\e$, which implies
\beq
\int_{\Om}D\psi \cdot F_A(Du^\e)  \geq 0 \text{\ \ \  for }\psi \leq 0,\  \psi \in W^{1,1}_0(\Om'^\e).
\eeq
For $\phi$ as in the beginning of the proof, choose  $\de>0$ small such that $\supp(\phi) \sub \Om'^\de$ and restrict $\e\leq \de$. Since $\Om'^\de \sub \Om'^\e$, by the elementary inequality \eqref{3.8} we have
\beq  \label{5.31}
0\, \leq\, \int_{\Om'^\de}F(Du^\e+D\psi) - \int_{\Om'^\de}F(Du^\e)
\eeq
for $\psi \in W^{1,1}_0(\Om'^\de)$,  $\psi \leq 0$ and all  $0<\e\leq \de$. Hence, $u^\e$ is a local subminimiser of \eqref{1.1} on $\Om'^\de$. Moreover, by Lemma \ref{le1} and \eqref{5.28}, we have $u^\e \overset{*}{\lharpoonup} u$ weakly$^*$ in $W^{1,\infty}(\Om'^\de)$ as $\e \ri 0$. By assumption \ref{1.2} and standard semicontinuity results (Dacorogna \cite{D}, p. 94), we have
\beq \label{5.32}
E\big(u,\Om'^\de\big)\, \leq\, \underset{\e\ri 0}{\lim\inf} \, E\big(u^\e,\Om'^\de \big).
\eeq
We now choose $\si>0$ small and define the cut-off function
\beq  \label{5.33}
\ze^\si\, :=\, \min \left\{\frac{1}{\si}\dist \big(\, \cdot\, , \p \Om'^\de \big),1\right\}.
\eeq
Then $\ze^\si \in W^{1,\infty}_0(\Om'^\de)$, $\ze^\si \equiv 1$ on the inner $\si$-neighborhood of $\Om'^\de$
\beq
\left\{x\in \Om'^\de \ \big| \ \dist(x,\p \Om'^\de)>\si \right\}
\eeq
 and $|D\ze^\si|\leq C/\si$ for a $C>0$ independent of $\si$. We select as $\psi$ the function
\beq \label{5.35}
\psi := \ze^\si\big(u-u^\e \big)\, +\, \phi,
\eeq
which is admissible since by Lemma \ref{le1} we have $u-u^\e\leq 0$.  We set
\beq \label{5.36a}
M\, :=\, 2\|Du\|_{L^\infty(\Om')}\, +\, \frac{C}{\si}\|u-u^\e\|_{C^0(\Om')}  \, +\, \|D\phi\|_{L^\infty(\Om')}.
\eeq
Then, since $F \in C^1(\R^n)$, by \eqref{5.31} we have
\begin{align} \label{5.37}
\int_{\Om'^\de}F(Du^\e) & \leq \int_{\Om'^\de}F\Big( (1-\ze^\si)Du^\e +\ze^\si Du + (u-u^\e)D\ze^\si+D\phi\Big) \nonumber\\
&\leq  \int_{\Om'^\de} F\big(\ze^\si Du +D\phi\big) \nonumber\\
&\ \ \ +\, \max_{\overline{\mB_M(0)}}|F_A|\int_{\Om'^\de}\big| (1-\ze^\si)Du^\e + (u-u^\e)D\ze^\si \big|\\
&\leq  \int_{\Om'^\de} F\big(\ze^\si Du +D\phi\big) \nonumber\\
&\ \ \ +\, \max_{\overline{\mB_M(0)}}|F_A|  
\left\{ \|Du\|_{L^\infty(\Om')} \int_{\Om'^\de} (1-\ze^\si) \, +\, \big|\Om'\big|  \frac{C}{\si}\|u-u^\e\|_{C^0(\Om')} \right\}. \nonumber
\end{align}
By \eqref{5.36a}, $M$ is independent of $\e$. Moreover, we have that $\ze^\si \nearrow 1$ a.e.\ on $\Om'^\de$ as $\si \ri 0$ and also the bound
\beq
F\big(\ze^\si Du +D\phi\big)\, \leq\, \max_{\overline{\mB_N(0)}}F,\ \ \ \ N:=\|Du\|_{L^\infty(\Om')} +\|D\phi\|_{L^\infty(\Om')} .
\eeq
By letting $\e \ri0$ and then $\si \ri 0$ in \eqref{5.37}, the Dominated Convergence theorem implies
\beq \label{5.39}
\underset{\e\ri 0}{\lim\inf} \int_{\Om'^\de} F(Du^\e)\, \leq \, \int_{\Om'^\de} F(Du+D\phi).
\eeq
By combining \eqref{5.39} with \eqref{5.32} and letting $\de\searrow 0$ we get
\beq \label{5.40}
\int_{\Om'} F(Du)\, \leq \, \int_{\Om'} F(Du+D\phi).
\eeq
for any $\phi \in C^\infty_c(\Om')$, $\Om'\Subset \Om$, $\phi \leq 0$. On the other hand, by following the same method for $u$ which is a Feeble Viscosity Supersolution and utilising the flat inf-convolution \eqref{4.2}, we deduce that \eqref{5.40} holds for $\phi \geq 0$ in $\phi \in C^\infty_c(\Om')$. By splitting a general $\phi$ to positive and negative part $\phi^+-\phi^-$, we have that \eqref{5.40}  holds for any $\phi  \in C^\infty_c(\Om')$. By taking variations in the standard way, \eqref{5.40} implies
\beq \label{5.41}
\int_{\Om'}F_A(Du) \cdot D\phi =0
\eeq
and since $|F_A(Du)| \in L^\infty(\Om')$, we have that \eqref{5.41} holds for all $\phi \in W^{1,1}_0(\Om')$. By inequality \eqref{3.8}, $u$ is a local minimiser as well and satisfies \eqref{1.3}. The proposition under assumption \eqref{1.6} follows.

\ms

\noi \textbf{2nd case under assumption \eqref{1.7}.} We begin with the next 

\bc \label{cl1} Let $F\in C^1(\R^n)$ satisfy \eqref{1.7} with $r>1$ as in \eqref{1.7} and assume that $v\in (W^{1,r}_{\text{loc}}\cap L^\infty_{\text{loc}})(\Om)$ is a weak solution of
\beq  \label{5.42}
\Div\big(\xi F_A(Dv) \big) \geq 0
\eeq
on $\Om\sub \R^n$, where $\xi \in \{1,-1\}$. Then, for any ball $\mB_R$ such that $\mB_{2R}\Subset \Om$, we have
\beq  \label{5.43}
\|Dv\|_{L^r(\mB_R)}\, \leq\, C\|v\|_{L^\infty(\mB_{2R})},
\eeq
where $C>0$ depends only $F,R$.
\ec

The proof of the claim is rather standard (see e.g.\ \cite{GT}), but instead of quoting general results with complicated assumptions, we prefer to give a short direct proof.

\BPC \ref{cl1}. Fix $\mB_{2R} \Subset \Om$ and $\ze \in C^\infty_c(\mB_{2R})$ with $\chi_{\mB_R}\leq \ze \leq \chi_{\mB_{2R}}$. 
\beq \label{5.44}
\psi \, :=\,  \ze^r \Big( \xi  v + \|v\|_{L^\infty(\mB_{2R})} \Big).
\eeq
Then, $\psi \in W^{1,r}_0(\mB_{2R})$, $\psi\geq 0$. Since $|F_A(Dv)| \in L^{\frac{r}{r-1}}_{\text{loc}}(\Om)$, by \eqref{5.42} and \eqref{5.44} we have
\begin{align}
-\int_{\mB_{2R}} \xi F_A(Dv) \cdot \Big[\xi  \ze^r Dv + r\Big( \xi v +\|v\|_{L^\infty(\mB_{2R})} \Big) \ze^{r-1} D\ze \Big]\geq 0.
\end{align}
Since $\xi^2=1$, in view of \eqref{1.7} and Young inequality, for $\e>0$ we obtain
\begin{align}
\int_{\mB_{2R}} \ze^r |Dv|^r & \leq C\|v\|_{L^\infty(\mB_{2R})} \int_{\mB_{2R}}  |Dv|^{r-1}\ze^{r-1}|D\ze| \nonumber\\
& \leq C\|v\|_{L^\infty(\mB_{2R})}  \left\{ \e \int_{\mB_{2R}}  \Big(|Dv|^{r-1} \ze^{r-1} \Big)^{\frac{r}{r-1}}  + \frac{1}{\e^{r-1}} \int_{\mB_{2R}} |D\ze|^r  \right\}.
\end{align}
For $\e:=1/\big(2C\|v\|_{L^\infty(\mB_{2R})}\big)$ we immediately get $\|\ze Dv\|_{L^r(\mB_{2R})}\, \leq\, C\|v\|_{L^\infty(\mB_{2R})}$ which implies \eqref{5.43}.   \qed

\bc \label{cl2} Suppose that $F\in C^1(\R^n)$ satisfies \eqref{1.7} with $r>1$ as in \eqref{1.7}. Assume that $(v^i)_1^\infty \sub (W^{1,r}_{\text{loc}}\cap L^\infty_{\text{loc}})(\Om)$ is a sequence of uniformly bounded in (the Fr\'echet topology of) $L^\infty_{\text{loc}}(\Om)$ weak solutions to
\beq  \label{5.47}
\Div\big(F_A(Dv) \big) \geq 0
\eeq
on $\Om\sub \R^n$. Then, if $v^{i+1}\leq v^i$ a.e.\ on $\Om$, the sequence $(v^i)_1^\infty$ is strongly precompact in (the Fr\'echet topology of) $W^{1,r}_{\text{loc}}(\Om)$. Moreover, if $v^i \ri v$ along a subsequence, then the limit solves \eqref{5.47} as well.
\ec

There is also a symmetric result for the compactness of increasing sequences of weak supersolutions, which we refrain from stating.

\BPC \ref{cl2}. The idea is taken from \cite{L}. Fix $\mB_{2R}$ and $\ze$ as in the proof of Claim \ref{cl1}. By \eqref{5.47} and our assumption, $(v^i)_1^\infty$ is weakly bounded in $W^{1,r}_{\text{loc}}(\Om)$ and hence there is $v$ such that $v^i \lharpoonup v$ weakly in $W^{1,r}_{\text{loc}}(\Om)$ along a subsequence as $i\ri \infty$. For each $i\in \N$, consider the integral
\beq  \label{5.48}
I^i\, :=\, \int_{\mB_{2R}} \Big(F_A(Dv^i) - F_A(Dv) \Big) \cdot D\big(\ze(v^i-v)\big).
\eeq
Since $\psi := \ze(v^i-v)$ is in $W^{1,r}_0(\mB_{2R})$ and also $\psi\geq 0$, by utilising that $v^i$ is a weak solution of \eqref{5.47} and that $|F_A(Dv^i)| \in L^{\frac{r}{r-1}}_{\text{loc}}(\Om)$, we have
\beq  \label{5.49}
I^i\, \leq \, \int_{\mB_{2R}}  F_A(Dv)  \cdot \big(D(\ze v^i)-D(\ze v)\big).
\eeq
Since $D(\ze v^i)  \rightharpoonup  D(\ze v)$ in $L^r(\mB_{2R})$ as $i\ri \infty$, by \eqref{5.49} we have that $\underset{i\ri \infty}{\lim\sup}\, I^i\leq 0$. By \eqref{1.7}, \eqref{5.48} and Young inequality we have
\begin{align} \label{5.50}
I^i\, &=\, \int_{\mB_{2R}} (v^i-v)  D\ze \cdot  \Big(F_A(Dv^i) - F_A(Dv) \Big) \nonumber\\
&\ \ \ \ \ \ +\int_{\mB_{2R}} \ze\Big(F_A(Dv^i) - F_A(Dv) \Big) \cdot D\big(v^i-v\big)\\
&\geq \,- \|D\ze \|_{L^\infty(\mB_{2R})} \|v^i-v \|_{L^r(\mB_{2R})}  \Big(\|Dv^i \|_{L^r(\mB_{2R})} + \| Dv\|_{L^r(\mB_{2R})}\Big)  \nonumber\\
& \ \ \ \ \ \ + C\int_{\mB_{2R}} \ze|Dv^i -Dv|^r. \nonumber
\end{align}
Hence, by \eqref{5.50} we get $\underset{i\ri \infty}{\lim\sup}\,  \|Dv^i-Dv \|_{L^r(\mB_{R})} = 0$, as desired.  \qed

\ms

We may now complete the proof of Proposition \ref{pr2} under assumption \eqref{1.7}. Suppose $u \in (C^0 \cap W^{1,1}_{\text{loc}})(\Om)$ is a Feeble Viscosity Subsolution of \eqref{1.4} on $\Om$. Fix $\Om'\Subset \Om$ and $\psi \in W^{1,r}_c(\Om')$ with $\psi\geq 0$ where $r>1$ is as in \eqref{1.7} and observe that $u\in C^0(\overline{\Om'})$. Consider the increasing function $\Phi:=\Phi^\infty \in C^0(0,\infty)$ with $\Phi(0^+)=0$ defined by \eqref{5.16} in Lemma \ref{le3}, where we take $R:=\infty$ since by assumption \eqref{1.7} the Laplacian of $F$ is bounded at infinity.  

Consider for $\e>0$  the flat sup-convolution $u^\e$ of $u$ restricted on $\Om'$, as given by \eqref{4.1}, where as $\Th \in C^2(0,\infty)$ we take the function given by Lemma \ref{le2} for the selected $\Phi$.  Then $u^\e$ satisfies all the properties of Lemmas \ref{le1}, \ref{le2}. Hence the flatness estimate \eqref{4.30} holds. Moreover, by Lemma \ref{le1}, $u^\e$ is a strong subsolution of
\beq
F_{AA}(Du^\e):D^2u^\e \, \geq\, 0,
\eeq
a.e.\ on $\Om'^\e \set \{Du^\e=0\}$. Since $u^\e$ is semiconvex, it satisfies the assumptions of Lemma  \ref{le3}. Hence, by Lemma  \ref{le3} $u^\e$ is a locally Lipschitz weak subsolution of $\Div\big(F_A(Du^\e) \big) \geq 0$ on $\Om'^\e$, and hence for $\e$ small such that $\supp(\psi)\sub \Om'^\e$, we have
\beq \label{5.52}
-\int_{\Om'}F_A(Du^\e) \cdot D\psi \,\geq\, 0.
\eeq
By Lemma \ref{le1} we have $u^\e \ri u$ in $C^0(\overline{\Om'})$ as $\e\ri 0$ and hence by Claims \ref{cl1}, \ref{cl2}, we have that  $u^\e \ri u$ in $W^{1,r}_{\text{loc}}(\Om')$  as $\e\ri 0$, along a sequence. By passing to the limit as $\e \ri 0$ in \eqref{5.52} and utilising \eqref{3.8}, we obtain that $u$ is a weak subolution of \eqref{1.12} in $W^{1,r}_{\text{loc}}(\Om)$ and also satisfies
\beq \label{5.53}
E(u,\Om') \leq E(u+\phi,\Om'),\ \ \Om'\Subset \Om,\ \ \phi \in W^{1,r}_0(\Om'),\ \phi\leq 0.
\eeq
Since $u$ is a Feeble Viscosity Supersolution as well, by arguing symmetrically for the flat inf-convolution, we obtain that $u$ is a weak supersolution as well, and \eqref{5.53} holds for $\phi\geq0$ as well. Hence, the proposition follows and so does Theorem \ref{th1}.        \qed

\ms

\ms

\ms

\noi \textbf{Acknowledgement.} I am indebted to J. Manfredi for his selfless share of expertise on the subject. I would also like to thank V. Julin for our scientific discussions. Special thanks are due to the referees of this paper for the careful reading of this manuscript, whose comments improved both the appearance and the content of the paper. In particular, I thank one of the referees for spotting an error in an earlier version of the manuscript on the construction of flat sup-convolutions.

\ms

\bibliographystyle{amsplain}

\begin{thebibliography}{10}
\bibitem[B]{B} J.M. Ball,  \emph{Some open problems in elasticity} In Geometry, Mechanics, and Dynamics, pages 3-59, Springer, New York, 2002.

\bibitem[BJ]{BJ} N. Barron, R. Jensen, \emph{ Minimizing the $L\infty$ norm of the gradient with an energy constraint}, Communications in PDE 30 (2005), 1741-1772.

\bibitem[CF]{CF} M. Colombo, A. Figalli, \emph{Regularity results for very degenerate elliptic equations}, J. Math\'ematiques Pures et Appliqu\'ees, Volume 101, Issue 1, January 2014, Pages 94 - 117.

\bibitem[CIL]{CIL} M. G. Crandall, H. Ishii, P.L. Lions, \emph{User's Guide to Viscosity Solutions of 2nd Order Partial Differential Equations}, Bulletin of the AMS, Vol. 27, Nr 1, 1-67, 1992.

\bibitem[D]{D} B. Dacorogna,  \emph{Direct Methods in the Calculus of Variations}, 2nd Edition, Volume 78, Applied Mathematical Sciences, Springer, 2008.

\bibitem[Da]{Da} J.M. Danskin, \emph{The theory of min-max with application}, SIAM Journal on Applied Mathematics, 14 (1966), 641 - 664.

\bibitem[EG]{EG} L.C. Evans, R. F. Gariepy, \emph{Measure theory and fine properties of functions}, Studies in Advanced Mathematics, CRC press, 1992.

\bibitem[GT]{GT} D. Gilbarg, N. Trudinger, \emph{Elliptic Partial Differential Equations of Second Order}, 
Classics in Mathematics, reprint of the 1998 edition, Springer.

\bibitem[I]{I} H. Ishii, \emph{On the Equivalence of Two Notions of Weak Solutions, Viscosity Solutions and Distribution Solutions}Funkcialaj Ekvacioj, 38, 101-120, 1995.

\bibitem[IS]{IS} H. Ishii, P. Souganidis, \emph{Generalized motion of noncompact hypersurfaces with velocity having arbitrary growth on the curvature tensor}, Tohoku Math Journal, 47, 227-250, 1995..

\bibitem[L]{L} P. Lindqvist, \emph{On the notion and properties of $p$-superharmonic functions}, J. Reine Angew. Math., 365, 67 - 79, 1986.

\bibitem[JJ]{JJ} V. Julin, P. Juutinen, \emph{A new proof for the equivalence of weak and viscosity solutions for the $p$-Laplace equation}, Communications in PDE, 37, No 5, 934-946, 2012.

\bibitem[JL]{JL} P. Juutinen, P. Lindqvist, \emph{Removability of a Level Set for Solutions of Quasilinear Equations}, Communications in PDE, Vol. 30, Issue 3, 305 - 321, 2005.

\bibitem[JLM]{JLM} P. Juutinen, P. Lindqvist, J.J. Manfedi, \emph{On the equivalence of viscosity solutions and weak solutions for a quasilinear equation}, SIAM J. Math. Anal., Vol. 33, 699-717, 2001. 

\bibitem[JLP]{JLP} P. Juutinen, T. Lukkari, M. Parviainen, \emph{Equivalence of viscosity and weak solutions for the $p(x)$-Laplacian} Annales de l'Institut Henri Poincare (C) Non Linear Analysis, Vol. 27, Issue 6, 1471-1487, 2010. 

\bibitem[S]{S} W. Schirotzek, \emph{Nonsmooth Analysis}, Universitext, Springer, 2007.

\bibitem[SS]{SS} S. Sheffield, C.K. Smart, \emph{Vector Valued Optimal Lipschitz Extensions}, Comm. Pure Appl. Math., Vol. 65, Issue 1, January 2012, 128 - 154.

\bibitem[SaV]{SaV} F. Santambrogio, V. Vespri, \emph{Continuity in two dimensions for a very degenerate elliptic equation}, Nonlinear Analysis TMA, Volume 73, Issue 12, 15 December 2010, Pages 3832 - 3841.

\bibitem[SV]{SV} R. Servadei, E. Valdinoci, \emph{Weak and viscosity solutions of the fractional Laplace equation}, to appear in Publ. Mat.
\end{thebibliography}

\end{document}